\documentclass[10pt]{amsart}
\title[Boundary $C^*$-algebras of $\tilde A_2$ groups]
{On the K-theory of boundary $C^*$-algebras of
$\widetilde A_2$ groups}
\author[]{Oliver King}
\author[]{Guyan Robertson}
\address{School of Mathematics and Statistics, University of Newcastle, NE1 7RU, England, U.K.}
\email{o.h.king@ncl.ac.uk}
\email{a.g.robertson@ncl.ac.uk}
\subjclass{Primary 46L80; secondary 58B34, 51E24, 20G25}
\keywords{affine building, boundary, operator algebra}

\usepackage{amsmath}
\usepackage{amscd}
\usepackage{amssymb}

\input{prepictex}
\input{pictex}
\input{postpictex}
\chardef\bslash=`\\ 

\makeatletter
\def\verbatim{\interlinepenalty\@M \@verbatim
  \leftskip\@totalleftmargin\advance\leftskip2pc
  \frenchspacing\@vobeyspaces \@xverbatim}
\makeatother
\newtheorem{theorem}{Theorem}[section]
\newtheorem{corollary}[theorem]{Corollary}
\newtheorem{lemma}[theorem]{Lemma}
\newtheorem{proposition}[theorem]{Proposition}

\theoremstyle{definition}
\newtheorem{example}[theorem]{Example}

\newtheorem{remark}[theorem]{Remark}

\newcommand{\cl}[1]{{\mathcal{#1}}}
\newcommand{\bb}[1]{{\mathbb{#1}}}
\newcommand{\fk}[1]{{\mathfrak{#1}}}

\newcounter{picture}

\DeclareMathOperator{\Tr}{Tr}
\newcommand{\e}{{\varepsilon}}

\newcommand{\PGL}{{\text{\rm{PGL}}}}

\newcommand{\GL}{{\text{\rm{GL}}}}

\newcommand{\vol}{{\text{\rm{vol}}}}
\newcommand{\id}{{\bf 1}}

\begin{document}

\begin{abstract}
Let $\Gamma$ be an $\widetilde A_2$ subgroup of $\PGL_3(\bb K)$, where $\bb K$ is a local field
with residue field of order $q$. The module of coinvariants $C(\bb P^2_{\bb K},\bb Z)_{\Gamma}$ is
shown to be finite, where $\bb P^2_{\bb K}$ is the projective plane over $\bb K$. If the group $\Gamma$ is of Tits type
and if $q \not\equiv 1 \pmod {3}$ then
the exact value of the order of the class $[\id]_{K_0}$ in the K-theory of the (full) crossed product $C^*$-algebra $C(\Omega)\rtimes\Gamma$ is determined, where $\Omega$ is the Furstenberg boundary of $\PGL_3(\bb K)$. For groups of Tits type, this verifies a conjecture of G. Robertson and T. Steger.
\end{abstract}

\maketitle

\section{Introduction}\label{intro}

Let $\bb K$ be a nonarchimedean local field with residue field $k$ of order $q$ and uniformizer $\pi$.
Denote by $\bb P^2_{\bb K}$ the set of one dimensional subspaces of the vector space $\bb K^3$, i.e. the set of points in the projective plane over $\bb K$. Then
$\bb P^2_{\bb K}$ is a compact totally disconnected space with the topology inherited from $\bb K$, and there is a continuous action of $G=\PGL_3(\bb K)$ on $\bb P^2_{\bb K}$. The group $G$ also acts on its Furstenberg boundary $\Omega$, which is the space of maximal flags $(0)<V_1<V_2<\bb K^3$.

The Bruhat-Tits building $\Delta$ of $G$ is a topologically contractible simplicial complex
of dimension $2$, which is a union of apartments. Each apartment in $\Delta$ is a flat subcomplex isomorphic to a tessellation of $\bb R^2$ by equilateral triangles:
that is, an affine Coxeter complex of type $\widetilde A_2$. For this reason, $\Delta$ is referred
to as an affine building of type $\widetilde A_2$.

\refstepcounter{picture}
\begin{figure}[htbp]
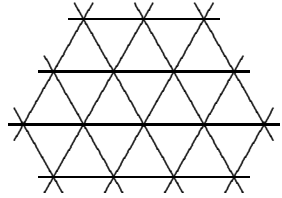
\label{A2tilde Coxeter}
\hfil
\centerline{
\beginpicture
\setcoordinatesystem units  <0.4cm, 0.7cm>        
\setplotarea x from -2.5 to 5, y from -1.5 to 2.5   
\putrule from -2.5   2     to  2.5  2
\putrule from  -3.5 1  to 3.5 1
\putrule from -4.5 0  to 4.5 0
\putrule from -3.5 -1  to  3.5 -1
\setlinear
\plot  -4.3 -0.3  -1.7 2.3 /
\plot -3.3 -1.3  0.3 2.3 /
\plot -1.3 -1.3  2.3 2.3 /
\plot  0.7 -1.3  3.3 1.3 /
\plot  2.7 -1.3  4.3 0.3 /
\plot  -4.3 0.3  -2.7 -1.3 /
\plot  -3.3 1.3   -0.7 -1.3 /
\plot  -2.3 2.3   1.3 -1.3 /
\plot  -0.3 2.3   3.3 -1.3 /
\plot  1.7 2.3   4.3 -0.3 /
\endpicture
}
\hfil
\caption{The Coxeter complex of type $\widetilde A_2$}
\end{figure}

\begin{figure}[htbp]
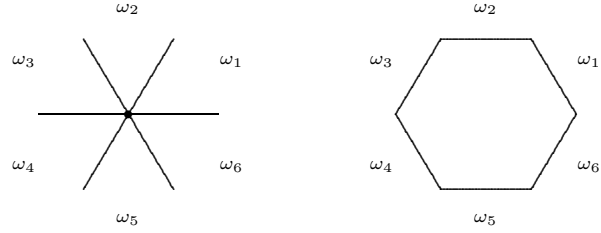
\label{A2 Coxeter}
\hfil
\centerline{
\beginpicture
\setcoordinatesystem units <0.6cm, 1cm>
\setplotarea  x from -3 to 3,  y from -1.5 to 1.4
\putrule from -2 0 to 2 0
\setlinear \plot  -1 -1  0 0  1 1  /
\setlinear \plot  1 -1  0 0  -1 1 /
\put {$_{\bullet}$} at 0 0
\put {$_{\omega_1}$} at 2.3 0.7
\put {$_{\omega_3}$} at -2.3 0.7
\put {$_{\omega_4}$} at -2.3 -0.7
\put {$_{\omega_6}$} at 2.3 -0.7
\put {$_{\omega_2}$} at 0   1.4
\put {$_{\omega_5}$} at 0 -1.4
\endpicture
\qquad \qquad
\beginpicture
\setcoordinatesystem units <0.6cm, 1cm>
\setplotarea  x from -2 to 3,  y from -1.5 to 1.4
\setlinear \plot   1 1  -1 1  -2 0  -1 -1   1 -1   2 0  1 1 /
\put {$_{\omega_1}$} at 2.3 0.7
\put {$_{\omega_3}$} at -2.3 0.7
\put {$_{\omega_4}$} at -2.3 -0.7
\put {$_{\omega_6}$} at 2.3 -0.7
\put {$_{\omega_2}$} at 0   1.4
\put {$_{\omega_5}$} at 0 -1.4
\endpicture
}
\hfil
\caption{The six boundary points of an apartment}
\end{figure}

An apartment in $\Delta$ is a union of six sectors, based at a fixed vertex
(Figure \ref{A2 Coxeter}). These six sectors correspond to six points $\omega_i$ in the Furstenberg boundary $\Omega$. They also represent the six edges of a hexagon in the spherical building at infinity $\Delta^\infty$. This hexagon (which is a spherical Coxeter complex of type $A_2$) is an apartment in $\Delta^\infty$. For this reason, $\Delta^\infty$ is referred to as a spherical building of type $A_2$.
More details are provided in Section \ref{boundary} below.

\refstepcounter{picture}

Let $\Gamma$ be a lattice subgroup of $G$. The abelian group $C(\bb P^2_{\bb K},\bb Z)$ of continuous integer-valued functions on $\bb P^2_{\bb K}$ is a $\Gamma$-module and the module of coinvariants $C(\bb P^2_{\bb K},\bb Z)_{\Gamma}= H_0(\Gamma; C(\bb P^2_{\bb K},\bb Z))$ is a finitely generated group.
Now suppose that $\Gamma$ acts freely and transitively on the vertex set of the Bruhat-Tits building of $G$, i.e. $\Gamma$ is an $\widetilde A_2$ group \cite{cmsz}.
Such a group is a natural analogue of a free group, which acts freely and transitively
on the vertex set of a tree (which is a building of type $\widetilde A_1$).
We prove that $C(\bb P^2_{\bb K},\bb Z)_{\Gamma}$ is a finite group and that the class $[\id]$ in $C(\bb P^2_{\bb K},\bb Z)_{\Gamma}$ has order
bounded by $q^2-1$.

An $\widetilde A_2$ group of \textit{Tits type} has a presentation with large automorphism group
\cite[II, Sections 4,5]{cmsz}.
If $\Gamma$ is of Tits type, and if $q \not\equiv 1 \pmod {3}$
 it is shown that the class $[\id]$ in $C(\bb P^2_{\bb K},\bb Z)_{\Gamma}$ has order exactly $q-1$.
For groups of Tits type, this proves a conjecture of \cite{rs3} regarding the
the order of the class $[\id]_{K_0}$ in the K-theory of boundary $C^*$-algebras.

\begin{theorem}\label{main}
   Let $\Gamma$ be an $\widetilde A_2$ group of Tits type and let $\fk A_{\Gamma}$ be the (full) crossed product $C^*$-algebra $C(\Omega)\rtimes\Gamma$.
  If $q \not\equiv 1 \pmod {3}$ then the order of $[\id]_{K_0}$ in $K_0(\fk A_\Gamma)$ is precisely $q-1$.
\end{theorem}

It is worth noting that the full crossed product
$C^*$-algebra $C(\Omega)\rtimes\Gamma$ coincides with the reduced crossed product, since the action of $\Gamma$
on $\Omega$ is amenable \cite[Section 4.2]{rs1}.
For $q \equiv 1 \pmod {3}$, the conjectured order of $[\id]_{K_0}$ is $\frac{q-1}{3}$,and this is verified for $q\le 31$. The element $[\id]_{K_0}$ is significant because the $C^*$-algebras $\fk A_{\Gamma}$ are classified up to isomorphism by their two $K$-groups, together with the class $[\id]_{K_0}$ \cite[Remark 6.5]{rs2}.
Many of the results are proved in a more general context, where the $\widetilde A_2$ group $\Gamma$ is not necessarily a lattice subgroup of $\PGL_3(\bb K)$.

\section{Background}
\subsection{The Bruhat-Tits building of $\PGL_3(\bb K)$}\label{building}

Given a local field $\bb K$, with discrete valuation $v~:~\bb K^\times\to\bb Z$, let
$\cl O=\{x\in \bb K:v(x)\ge0\}$ and let $\pi\in \bb K$ satisfy $v(\pi)=1$.
A \textit{lattice} $L$ is a rank-3 $\cl O$-submodule of $\bb K^3$.
In other words $L=\cl O v_1 + \cl O v_2 + \cl O v_3$, for some basis
$\{v_1, v_2, v_3\}$ of $\bb K^3$.
Two lattices $L_1$ and $L_2$ are \textit{equivalent} if $L_1=uL_2$ for some $u\in \bb K^\times$.
The Bruhat-Tits building $\Delta$ of $\PGL_3(\bb K)$ is a two dimensional simplicial complex whose vertices are equivalence classes of lattices in $\bb K^3$.
Two lattice classes $[L_0], [L_1]$ are \textit{adjacent} if, for suitable representatives
 $L_0, L_1$, we have $L_0 \subset L_1 \subset \pi^{-1} L_0$. A \textit{simplex} is a set of pairwise adjacent lattice classes. The maximal simplices (chambers) are  the sets $\{[L_0], [L_1],[L_2]\}$ where
$L_0 \subset L_1\subset L_2 \subset \pi^{-1} L_0$. These inclusions determine a canonical ordering of the vertices in a chamber, up to cyclic permutation.

Each vertex $v$ of $\Delta$ has a {\it type} $\tau(v) \in \bb Z/3\bb Z$, and each chamber of $\Delta$ has exactly one vertex of each type.
If the Haar measure on $\bb K^3$ is normalized so that $\cl O^3$ has measure $1$ then the type map may be defined by
$\tau([L])=\vol(L) + 3\bb Z$. The cyclic ordering of the vertices of a chamber coincides with the natural ordering given by the vertex types.
The building $\Delta$ is of type $\widetilde A_2$ and there is a natural action of $\PGL_3(\bb K )$ on $\Delta$
induced by the action of $\GL_3(\bb K)$  on the set of lattices.

\subsection{$\widetilde A_2$ groups}\label{A2tildegroups}

More generally, consider any locally finite building $\Delta$ of type $\widetilde A_2$.  Each vertex $v$ of $\Delta$ is labeled with a type $\tau (v) \in \bb Z/3\bb Z$,
and each chamber has exactly one vertex of each type.
Each edge $e$ is directed, with initial vertex $o(e)$ of type $i$ and final vertex $t(e)$ of type $i+1$.
An automorphism $\alpha$ of $\Delta$ is said to be \textit{type rotating} if there exists $i \in
\bb Z/3\bb Z$ such that $\tau(\alpha(v)) = \tau(v)+i$ for all vertices $v \in \Delta$.

Suppose that $\Gamma$ is a group of type rotating automorphisms of $\Delta$, which acts freely and transitively on the vertex set of $\Delta$. Such a group $\Gamma$ is called an $\widetilde A_2$ group. The theory of $\widetilde A_2$ groups has been developed in \cite{cmsz} and some, but not all, $\widetilde A_2$ groups embed as lattice subgroups of $\PGL_3(\bb K)$.
Any $\widetilde A_2$ group can be constructed as follows \cite [I, Section 3]{cmsz}.
Let $(P,L)$ be a projective plane of order $q$. There are $q^2 + q + 1$
points
 (elements of $P$) and $q^2+q+1$
lines (elements of $L$).
Let $\lambda : P \rightarrow L$ be a bijection.
A {\it triangle presentation}  compatible with $\lambda$ is a set ${\mathcal T}$ of ordered triples $(x, y, z)$
where $x, y, z \in P$, with the following properties.
\begin{itemize}
\item[(i)]  Given $x, y \in P$, then $(x, y, z) \in {\mathcal T}$ for some
$z \in P$ if and only if $y$ and $\lambda(x)$ are incident, i.e. $y \in
\lambda(x)$.

\item[(ii)] $(x, y, z) \in {\mathcal T} \Rightarrow (y, z, x) \in {\mathcal T}$.

\item[(iii)]  Given $x, y \in P$, then $(x, y, z) \in {\mathcal T}$ for at
most one $z \in P$.
\end{itemize}
In \cite{cmsz} there is displayed a complete list of triangle
presentations for $q = 2$
and $q = 3$.
Given a triangle presentation $\cl T$, one can form the group
\begin{equation}\label{rel}
\Gamma=\Gamma_{\cl T} = \big\langle P \ |\  xyz = 1 \hbox { for } (x,
y, z) \in {\mathcal T}
\big \rangle .
\end{equation}
\noindent The Cayley graph of $\Gamma$ with respect to the
generating set $P$ is the $1$-skeleton of a building $\Delta$ of type $\widetilde A_2$.
Vertices are elements of $\Gamma$ and a directed edge of the form $e=(a,ax)$ with $a\in\Gamma$ is labeled by a generator $x\in P$.
Denote by $E_+$ the set of directed edges of $\Delta$.

The link of a vertex $a$ of $\Delta$ is the incidence graph of the projective plane $(P,L)$, where
the lines in $L$ correspond to the inverses in $\Gamma$ of the
generators in $P$. In other words,
$\lambda(x) = x^{-1}$ for $x \in P$.

\refstepcounter{picture}
\begin{figure}[htbp]
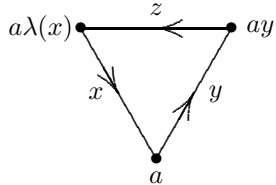
\label{A4}
\hfil
\centerline{
\beginpicture
\setcoordinatesystem units <1cm, 1.732cm>
\put {$\bullet$} at 0 0
\multiput {$\bullet$} at -1 1 *1 2 0 /
\put {$a$} [t] at 0 -0.1
\put {$a\lambda(x)$} [r] at -1.1 1
\put {$ay$} [l] at 1.2 1
\arrow <10pt> [.2, .67] from  0.2 1 to 0 1
\arrow <10pt> [.2, .67] from  -0.7 0.7 to -0.5 0.5
\arrow <10pt> [.2, .67] from  0.3 0.3 to 0.5 0.5
\put {$x$} [r ] at -0.7 0.5
\put {$y$} [ l] at 0.7 0.5
\put {$z$} [ b] at 0 1.1
\putrule from -1 1 to 1 1
\setlinear \plot -1 1 0 0 1 1 /
\endpicture
}
\hfil
\caption{A chamber based at a vertex $a$}
\end{figure}

\section{The group $A_{\cl T}$}

The main results will be proved by defining a homomorphism from an abelian group
$\cl A_{\cl T}$ onto the module of coinvariants.
Suppose that $\cl T$ is a triangle presentation with associated $\widetilde A_2$ group $\Gamma=\Gamma_{\cl T}$. Define $A_{\cl T}$ to be the abelian group generated by $P\cup\{\varepsilon\}$ subject to the relations

\begin{subequations}\label{rel1}
\begin{eqnarray}
\sum_{y\notin \lambda(x)}y &=& x, \quad\, x\in P; \label{a}\\
x+y+z &=& \varepsilon, \quad (x,y,z)\in \cl T; \label{b}\\
\sum_{x\in P}x &=& \varepsilon.\label{c}
\end{eqnarray}
\end{subequations}
It follows immediately from (\ref{a}) and (\ref{c}) that, for each $x\in P$,
$$
\e=\sum_{y\in P}y=\sum_{y\notin \lambda(x)}y +  \sum_{y\in \lambda(x)}y
= x + \sum_{y\in \lambda(x)}y.
$$
If we define $\widehat x = \displaystyle\sum_{y\in \lambda(x)}y$ then we obtain
\begin{equation}\label{d}
 x + \widehat x = \e, \quad\, x\in P.
\end{equation}
The group $A_{\cl T}$ has an alternative presentation  with relations (\ref{b}), (\ref{c}), (\ref{d}).
 It turns out that $A_{\cl T}$ is finite. First of all, observe that
the element $\e$ has finite order.

\begin{lemma}\label{q2}
In the group $A_{\cl T}$,  $(q^2-1)\e = 0$.
\end{lemma}

\begin{proof}
  Define a $\{0,1\}$-matrix $M$ by $M(x,y)=1 \Leftrightarrow y\notin \lambda(x)$.
Then
  \begin{equation*}
\begin{split}
  \e&=\sum_{x\in P}x=\sum_{x\in P}\sum_{y\in P} M(x,y)y \\
   &= \sum_{y\in P}(\sum_{x\in P}M(x,y))y = \sum_{y\in P}q^2 y \\
   &= q^2\e
\end{split}
\end{equation*}
\end{proof}

\begin{proposition}\label{finite}
  The group $A_{\cl T}$ is finite.
\end{proposition}

\begin{proof}
  The $\widetilde A_2$ group $\Gamma$ has Kazhdan's property T \cite{cms,Zuk}. It follows that the abelianization $\Gamma_{ab}=H_1(\Gamma, \bb Z)$ is finite \cite[Corollary 1.3.6]{bhv}.
  By relation (\ref{b}), $A_{\cl T}/\langle \e \rangle$ is a quotient of $\Gamma_{ab}$. Thus  $A_{\cl T}/\langle \e \rangle$ is finite.
Therefore so also is $A_{\cl T}$, since $\langle \e \rangle$ is finite, by Lemma \ref{q2}.
\end{proof}

\bigskip

\noindent\textbf{Question}. What is the order of $\e$ in $A_{\cl T}$?

\bigskip

\begin{lemma}\label{lower bound}
The order of $\e$ is at least $\frac{q-1}{(q-1,3)}$.
\end{lemma}

\begin{proof}
The map $f: P\cup \{\e\} \to \bb Z_{q^2-1}$ defined by $f(\e)=3(q+1)$ and
$f(x)=q+1$ for each $x\in P$, extends to a homomorphism $f: A_{\cl T} \to \bb Z_{q^2-1}$.
For $f$ preserves the relations (\ref{rel1}), since
$$\sum_{y\notin\lambda(x)}f(y) = q^2(q+1)=q+1=f(x), \quad\, x\in P,$$
$$f(x)+f(y)+f(z) = 3(q+1) =f(\e), \quad (x,y,z)\in \cl T,$$
and
$$\sum_{x\in P}f(x) = (q^2+q+1)(q+1)=3(q+1)=f(\e).$$
Since $f(\e)=3(q+1)$, the order of $f(\e)$ is
$$\frac{q^2-1}{(q^2-1,3(q+1))} = \frac{q^2-1}{(q+1)(q-1,3)}= \frac{q-1}{(q-1,3)}.$$
Therefore $\e$ has order at least $\frac{q-1}{(q-1,3)}$.
\end{proof}

\bigskip

\begin{proposition}\label{improvement}
Suppose that there exists a subset  $\cl M \subset \cl T$
such that each element in $P$ occurs precisely 3 times in the triples belonging to $\cl M$.
Then $(q-1)\e =0$.
\end{proposition}

\begin{proof}
Since $\# P = q^2+q+1$, we have $\# \cl M = q^2+q+1$. Therefore
\begin{equation*}
\begin{split}
  3\e & =  3\displaystyle \sum_{x\in P} x
   = \sum_{(x,y,z)\in \cl M} (x+y+z) \\
   &= \sum_{(x,y,z)\in \cl M} \e =  (q^2+q+1)\e.
\end{split}
\end{equation*}
The result follows, since $q^2\e=\e$,  by Lemma \ref{q2}
\end{proof}

\begin{remark}
Such a set $\cl M$ exists if $\cl T$ is a triangle presentation of any torsion free $\widetilde A_2$ group with $q=2$. Thus $\e=0$ if $q=2$. To see this, observe that for each $x_0\in P$ there are $3$ elements $y\in P$ such that $y\in \lambda(x_0)$.
Thus there are $3$ elements in $\cl T$ of the form $(x_0,y,z)$, where $z$ is uniquely determined by $x_0$ and $y$.
The $7$ possible $x_0\in P$ give rise to $7.3=21$ elements of $\cl T$.
Thus $\#\cl T=21$, and each element of $P$ occurs $3+3+3=9$ times in the triples belonging to $\cl T$. Since $\cl T$ contains no triple of the form $(x,x,x)$, the orbit of an element of $\cl T$ under cyclic permutations contains 3 elements.
Choosing one element of $\cl T$ from each such orbit, we obtain a set $\cl M$ containing $7$ elements of $\cl T$, in which each generator appears $3$ times.
\end{remark}
A weak form of Proposition \ref{improvement} is the following.

\begin{corollary}\label{improvement3}
Suppose that there exists a subset  $\cl M \subset \cl T$
such that each of the three  coordinate projections from $\cl M$ onto $P$ is bijective. Then
 $$(q-1)\e=0.$$
\end{corollary}
The rest of this section provides examples where Corollary \ref {improvement3} applies.

\subsection{Invariant triangle presentations}

\quad Consider the Desarguesian projective plane $(P,L)=PG(2,q)$. The points and lines are the 1- and 2-dimensional subspaces, respectively, of a 3-dimensional vector space $V$ over $\bb F_q$, with
incidence being inclusion.
The field $\bb F_q$ is a subfield of $\bb F_{q^3}$ and so
$\bb F_{q^3}$ is a 3-dimensional vector space over $\bb F_q$.
We may therefore identify $P$ with $\bb F_{q^3}^\times/\bb F_q^\times$.
Now $\bb F_{q^3}^\times$ is a cyclic group, with generator $\zeta$, say,
and $\bb F_q^\times$ is the subgroup generated by $\zeta^{1+q+q^2}$
(the unique subgroup with $q-1$ elements).
Thus $P$ is  a cyclic group of order $1+q+q^2$ generated by the element
$\sigma=\bb F_q^\times\zeta$.
Multiplication by $\zeta$ is an $\bb F_q$-linear transformation and so
multiplication by $\sigma$ is a collineation of $(P,L)$.
Thus $\sigma$ generates a cyclic collineation group of order $q^2+q+1$ which acts regularly on $P$ (and on $L$). A group with these properties is called a {Singer group}, after \cite{sin}.

\begin{remark}
 It is not known whether a projective plane with a Singer group is necessarily Desarguesian.
 Also, there are no known $\widetilde A_2$ groups for which the underlying projective
   plane $(P,L)$ is not Desarguesian.
\end{remark}

Call a triangle presentation S-invariant if it is invariant under a Singer group
of collineations of $(P,L)$.
Such presentations are studied in \cite[I, Section 5]{cmsz}.
It is shown in \cite[I, Theorem 5.1]{cmsz} that there are exactly $2^{(q+1)/3}$,
$2^{q/3}$ or $3(2^{(q-1)/3})$ distinct S-invariant triangle presentations, according as $q\equiv-1$, $0$
or~$1$ mod~$3$, respectively.
These presentations are constructed explicitly in  \cite{cmsz} and all of them are also invariant under the
Frobenius collineation $x\mapsto x^q$.
It is not known whether all the corresponding $\widetilde A_2$ groups embed as lattices
in $\PGL_3(\bb K)$.

\begin{remark}\label{inv}
If $\lambda : P \to L$ is a bijection and $\cl T$ is an S-invariant triangle presentation compatible
with $\lambda$ then $y\in \lambda(x) \leftrightarrow \sigma y\in \lambda(\sigma x)$. Therefore
the presentation (\ref{rel1}) of the group $\cl A_\cl T$ is invariant under the Singer group
$<\sigma>$. It is also  invariant under the collineation $x\mapsto x^q$.
\end{remark}

\begin{corollary}\label{S-invariant}
  If the triangle presentation $\cl T$ is S-invariant then, in the group $A_{\cl T}$,
  $(q-1)\e=0$.
\end{corollary}

\begin{proof}
 Fix a triple $(x_0,y_0,z_0)\in \cl T$. The orbit $\cl M_0$ of $(x_0,y_0,z_0)$ under the Singer group satisfies the hypotheses of Corollary \ref{improvement3}.
\end{proof}

\begin{remark}\label{conjecture}
It follows from Lemma \ref{lower bound} that if $\cl T$ is S-invariant and $q \not\equiv 1 \pmod {3}$ then the order of $\e$ is $q-1$.
Computations on the family of triangle presentations described below
confirm that, up to $q=32$, the order of $\e$ is $\frac{q-1}{(q-1,3)}$.
We conjecture that this is true for all S-invariant triangle presentations and for all $q$.
\end{remark}

\subsection{S-invariant triangle presentations of Tits type \cite[I Section 4]{cmsz}}\label{appendix}

The mapping  $x\mapsto x^q$ is an automorphism of $\bb F_{q^3}$ over $\bb F_q$, i.e. it has fixed field $\bb F_q$.
The trace $\Tr : \bb F_{q^3} \to \bb F_q$ is the $\bb F_q$-linear mapping defined by
$\Tr(a)=a + a^q + a^{q^2}$.
Now
$(x_0,y_0)\mapsto\Tr(x_0y_0)$ is a regular symmetric bilinear form on $\bb F_{q^3}$, and so
the map
$$V\mapsto V^\perp =\{y_0\in \bb F_{q^3}:\Tr(x_0y_0)=0
\text{  for all  } x_0\in V\}$$
is a bijection $P\to L$.

Let  $x={\bb F}_q^\times x_0\in P=\bb F_{q^3}^\times/\bb F_q^\times$.
Write $\Tr(x)=0$ if $\Tr(x_0)=0$ and write $x^\bot$ for
$({\bb F}_qx_0)^\bot$. We identify any line, i.e., 2-dimensional
subspace, with the set of its 1-dimensional subspaces. Under these
identifications, the lines are the subsets $x^\bot=\{y\in P:\Tr(xy)=0\}$
of the group~$P$.
For $x\in P$, let
$$
\lambda_0(x)=\left(\frac{1}{x}\right)^{\bot}=\left\{y\in P:\Tr\left(\frac{y}{x}\right)=0\right\}.
$$
This defines a point-line correspondence $\lambda_0:P\to L$. The
following set of triples is a triangle presentation $\cl T_0$ compatible with
$\lambda_0$:
\begin{equation}\label{regularpresentation}
\cl T_0=\{(x,x\xi,x\xi^{q+1})~:~x,\xi\in P {\rm\ and\ } \Tr(\xi)=0\}.
\end{equation}
The only nontrivial thing to check is that $\cl T_0$ is invariant under cyclic permutations.
If $(x,x\xi,x\xi^{q+1})\in \cl T_0$ then
$$
(x\xi,x\xi^{q+1},x)=(x\xi,(x\xi)\xi^q,(x\xi)(\xi^q)^{q+1}) \in \cl T_0
$$
since $\Tr(\xi_0^q)=\Tr(\xi_0)$ and
$\xi_0^{1+q+q^2}\in{\bb F}_q^\times$ for each $\xi_0\in{\bb F}_{q^3}^\times$.

The group $\Gamma_{\cl T}$ associated with the presentation (\ref{regularpresentation})
is called a group of Tits type \cite[I  Section 4]{cmsz}.
The presentation (\ref{regularpresentation}) is clearly invariant under multiplication by
any $x\in P$. That is, it is S-invariant. It is also easy to see that (\ref{regularpresentation})
is invariant under the Frobenius collineation $x\mapsto x^q$.

\begin{remark}\label{T0'}
Consider the related triangle presentation
\begin{equation}\label{semiregularpresentation}
\begin{split}
\cl T_0'=&\left\{\left(\frac{1}{z},\frac{1}{y},\frac{1}{x}\right)~:~(x,y,z)\in\cl T_0\right\}\\
=&\{(x,x\xi,x\xi^{q^2+1})~:~x,\xi\in P {\rm\ and\ } \Tr(\xi)=0\}
\end{split}
\end{equation}
which is compatible with $\lambda_0$ and S-invariant.
For $q=2,3$, the two presentations $\cl T_0, \cl T_0'$ are the only S-invariant triangle presentations
\cite[II,p.165,Remark]{cmsz}.

The group $\Gamma_{\cl T_0'}$ is isomorphic to $\Gamma_{\cl T_0}$, via the map
 $x\mapsto \left(\frac{1}{x}\right)$,
but the groups $A_{\cl T_0'}$, $A_{\cl T_0}$ are not isomorphic. In fact, writing $q=p^r$, $p$ prime, the MAGMA computer algebra package calculates the following expressions for the groups $A_{\cl T_0}$, $A_{\cl T_0'}$ for $2\le q\le 32$.

\bigskip

\centerline{
{\small
\begin{tabular}{|l|l|l|}
\hline
                                     & $q \not\equiv 1 \pmod {3}$ & $q \equiv 1 \pmod {3}$    \\
\hline
 $A_{\cl T_0}$          & $\bb Z_{q-1}$ & $\bb Z_{q-1}\oplus\bb Z_3$    \\
$A_{\cl T_0'}$           & $(\bb Z_p)^{3r}\oplus \bb Z_{q-1}$ & $(\bb Z_p)^{3r}\oplus \bb Z_{q-1}\oplus \bb Z_3$    \\ \hline
\end{tabular}
}
}
\bigskip
\noindent In all cases, the order of $\e$ is precisely $\frac{q-1}{(q-1,3)}$.
\end{remark}

\begin{example}
We verify the stated value of $\e$ in the case $q=4$.
The set $P$ is identified with $\bb F_{64}^\times/\bb F_4^\times$.
Write $\bb F_4=\bb F_2(\alpha)$, where $\alpha^2+\alpha+1=0$ and
$\bb F_{64}=\bb F_2(\alpha)(\delta)$, where $\delta^3+\delta^2+\alpha\delta+\alpha+1=0$.
Then $\delta^{21}=\alpha+1$ and $\delta^{42}=\alpha$ are the roots of $x^2+x+1$.
A short calculation shows that $(\delta^7)^4~=~(\alpha+1)\delta^7$, so
that the fixed points of the Frobenius collineation $x\mapsto x^4$ on $P=\bb F_{64}^\times/\bb F_4^\times$ are
$\dot{1}, \omega, \omega^2$, where $\dot{1}=\bb F_4^\times$ and  $\omega=\delta^7\bb F_4^\times$.
It is easy to check that $\{\xi\in P : \Tr(\xi)=0\}=\{\omega, \omega^2, s, s^4, s^{16}\}$,
where $s=\delta^3\bb F_4^\times$.
It follows from relation (\ref{d}) in the group $A_{\cl T_0}$, with $x=\dot{1}$, that
\begin{equation*}
\begin{split}
\e&= \dot{1}+(\omega+ \omega^2 + s+ s^4+ s^{16}) \\
&= (\dot{1}+\omega+ \omega^2) + (s+ s^4+ s^{16}) \\
&= \e + (s+ s^4+ s^{16}).
\end{split}
\end{equation*}
In the last line, the fact that $\dot{1}+\omega+ \omega^2=\e$ follows from (\ref{b}), since $(\dot{1},\omega,\omega^2)\in \cl T_0$ by (\ref{regularpresentation}). Therefore
\begin{equation*}
\begin{split}
\e+s^3&= \e + (s^3+ s^4+ s) + s^{16} \\
&= \e + (s^3+ s^4+ s^8) + s^2 \\
&= \e + \e + s^2,
\end{split}
\end{equation*}
using the facts that $s^8=\delta^{24}\bb F_4^\times=s$ and $(s^3,s^4,s^8)\in \cl T_0$. Hence
$s^3=\e+s^2$. By invariance of the presentation $\cl T_0$ (Remark \ref{inv}), we deduce that
$s-\dot{1}=\e$ and therefore that $\e=s^{16}-\dot{1}=s^2-\dot{1}$.
Thus $s^2=s$ and, by invariance, $s=\dot{1}$, which gives $\e=0$.
It is also easy to see that $A_{\cl T_0}=\bb Z_3\oplus\bb Z_3$, with generators $\delta$ and $\delta^2$.
\end{example}

\subsection{More examples.}\label{not}
  Consider any triangle presentation $\cl T$.  If $\phi$ is an order-3 collineation such that $\cl T$ is fixed by the map
$x\mapsto \phi(x)$, then one obtains a
new triangle presentation
\begin{equation}\label{newrel}
\cl T^\phi = \big\{(x,\phi(y),\phi^2(z)) :  (x,y,z) \in \mathcal T
\big \}
\end{equation}
relative to the point--line correspondence $\phi\circ\lambda : P \rightarrow L$.
The corresponding group
\begin{equation}\label{newgp}
\Gamma_{\cl T^\phi} = \big\langle P \ |\  x\phi(y)\phi^2(z) = 1 \hbox { for } (x,
y, z) \in {\mathcal T}
\big \rangle
\end{equation}
is \textit{not} in general isomorphic to $\Gamma_{\cl T}$ \cite[II]{cmsz}.

One possible choice for $\phi$ is the collineation $x\mapsto x^q$.  Applied
to the triangle presentations (\ref{regularpresentation}), $\phi$ and $\phi^2$ give two more triangle presentations which are \textit{not} in general S-invariant.
If $q \equiv 1 \pmod {3}$ then another possible choice of $\phi$ is $x\mapsto \omega x$, where
$\omega=\sigma^{(q^2+q+1)/3}$.

\begin{corollary}\label{almost S-invariant}
  If the triangle presentation $\cl T$ is S-invariant and if
  $\phi$ is an order-3 collineation of $(P,L)$ such that $\cl T$ is fixed by the map
$x\mapsto \phi(x)$ then $(q-1)\e=0$ in $A_{\cl T^\phi}$.
\end{corollary}

\begin{proof}
Let $\cl M= \{(x,\phi(y),\phi^2(z)) : (x,y,z) \in {\cl M_0}\}$,
where $\cl M_0$ is defined in the proof of Corollary \ref {S-invariant}.
This satisfies the hypotheses of Corollary \ref{improvement3}.
\end{proof}

\section{The boundary action}\label{boundary}

 Associated with the building $\Delta$ is the building at infinity $\Delta^\infty$. This is a spherical building of type $A_2$ \cite[Theorem 8.24]{Wei}.
In the geometrical realization of $\Delta^\infty$, a point $\xi \in \Delta^\infty$ is an equivalence class of rays (subsets of $\Delta$ isometric to $[0,\infty)$),
  where two rays are equivalent if they are \textit{parallel}, i.e. at finite Hausdorff distance from each other \cite[11.8]{ab}.
For each vertex $x$ of $\Delta$ and each $\xi\in\Delta^\infty$, there is a unique ray
$[x,\xi)$ with initial vertex $x$ in the parallelism class of $\xi$.
A \emph{sector} is a $\frac{\pi}{3}$-angled sector made up of chambers in some apartment.
Each sector in  $\Delta$ determines a 1-simplex (chamber) of $\Delta^\infty$ whose points
are equivalence classes of rays in the sector. Two sectors determine the same chamber of $\Delta^\infty$
if and only if they contain a common subsector.

If $[x,\xi)$ is a sector wall then $\xi$ is a \textit{vertex} of $\Delta^\infty$. If the initial edge of this sector wall is $[x,y]$, and if $\tau(y) = \tau(x)+i$ then the vertex $\xi$ is said to be of type $i-1$.
This definition is independent of the base vertex $x$, since two parallel rays lie in a common apartment.
Denote by $\cl P$ the set of vertices $\xi\in\Delta^\infty$ of type $0$ and denote by
$\cl L$ the set of vertices $\eta\in\Delta^\infty$ of type $1$.
Then $({\cl P},\cl L)$ is a projective plane.
A point $\xi \in \cl P$ and a line $\eta \in \cl L$ are incident if and only if they are the two vertices
of a common chamber in  $\Delta^\infty$. This is equivalent to saying that there is a sector
with base vertex $x$ whose walls are the rays $[x,\xi)$ and $[x,\eta)$.

\refstepcounter{picture}
\begin{figure}[htbp]
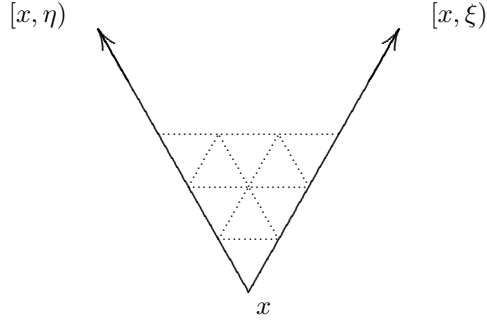
\label{A2}
\centerline{
\beginpicture
\setcoordinatesystem units <0.4cm,0.7cm>   
\setplotarea x from -4 to 4, y from  0 to  5  
\setlinear \plot -5 5  0 0  5 5  /
\arrow <8pt> [.2, .67] from  -4 4 to -5 5
\arrow <8pt> [.2, .67] from  4 4 to 5 5
\put {$[x,\xi)$} [ b] at 7 5
\put {$[x,\eta)$} [ b] at -7 5
\put {$x$} at 0.5 -0.3
\setdots <2.0pt>
\plot -1 1  0 2  1 1 /
\plot  -2 2  -1 3  0 2  1 3  2 2  /
\plot  -3 3    3 3 /
\plot  -1 1   1 1 /
\plot  -2 2   2 2 /
\endpicture
}
\caption{Sector walls}
\end{figure}

If $e\in E_+$, let $\Omega(e)$ denote the set of points $\xi\in {\cl P}$ which have representative rays with initial edge $e$.
That is,
\begin{equation*}
\Omega(e)=\{\xi\in {\cl P}:e\subset[o(e),\xi)\}\,.
\end{equation*}
The sets $\Omega(e)$, $e\in E_+$, form a basis of clopen sets for a
totally disconnected compact topology on $\cl P$.
The topological space $\cl P$ is called the \textit{minimal boundary} of $\Delta$.

\begin{remark}
If $\Delta$ is the Bruhat-Tits building of $\PGL_3(\bb K)$ then there is a natural identification
of $\cl P$ and $\bb P^2_{\bb K}$ as topological $G$-spaces.
If $\lambda \in \bb P^2_{\bb K}$ and $L_0$ is a lattice, define a sequence of lattices
inductively by $L_{i+1}= L_i +(\lambda \cap \pi^{-1}L_i)$. Then $L_i\subset L_{i+1}\subset\pi^{-1} L_i$ and
$\tau([L_{i+1}])=\tau([L_i])+1$, since $L_i$ is maximal in $L_{i+1}$. The sequence of vertices $[L_0], [L_1], [L_2], \dots$ defines a ray whose parallelism class $\xi_\lambda$ is an element of $\cl P$ and the map
$\lambda\mapsto\xi_\lambda$ is a bijection from $\bb P^2_{\bb K}$ onto $\cl P$.
If $e=([L_0],[L_1])$, where $L_0\subset L_1\subset \pi^{-1} L_0$ and $\tau([L_1])=\tau([L_0])+1$, then $\Omega(e)$
 may be identified with the set of lines $\lambda\in \bb P^2_{\bb K}$ such that $L_1=L_0 + (\lambda \cap \pi^{-1}L_0)$.
\end{remark}

\noindent If $v$ is a fixed vertex of $\Delta$, then ${\cl P}$ may be expressed as a disjoint union
\begin{equation}\label{boundary C}
{\cl P}=\bigsqcup_{o(e)=v}\Omega(e)\,.
\end{equation}
Also, if $e\in E_+$,
then $\Omega(e)$ can be expressed as a disjoint union
\begin{equation}\label{boundary A}
\Omega(e)=\bigsqcup_{\substack{o(e')=t(e)\\ \Omega(e') \subset \Omega(e)} }\Omega(e')\,.
\end{equation}

If $\Gamma$ is an $\widetilde A_2$ group acting on $\Delta$, then $\Gamma$ acts on ${\cl P}$, and the abelian group $C(\cl P,\bb Z)$ has the structure of a $\Gamma$-module, with
$(g\cdot f)(\xi)=f(g^{-1}\xi),\, g\in \Gamma,\, \xi\in\cl P$.  The module of $\Gamma$-coinvariants, $C(\cl P,\bb Z)_\Gamma$, is the quotient of
$C(\cl P,\bb Z)$ by the submodule generated by
$\{g\cdot f-f : g\in\Gamma, f\in C(\cl P,\bb Z)\}$.
If $f\in C(\cl P,\bb Z)$ then $[f]$ denotes its class in $C(\cl P,\bb Z)_\Gamma$.
Also, $\id$ denotes the constant function defined by $\id(\xi)=1$ for $\xi\in\cl P$.

If $e\in E_+$, let $\chi_e$ be the characteristic function of $\Omega(e)$.
For each $g\in\Gamma$, the functions $\chi_e$ and $g\cdot\chi_e=\chi_{ge}$ project to the same element in
$C(\cl P,\bb Z)_\Gamma$. Therefore, for any edge $e= (a,ax)$ with $a\in\Gamma$ and $x\in P$,
it makes sense to denote by $[x]$ the class of $\chi_e$ in $C(\cl P,\bb Z)_\Gamma$.

   Suppose that $e,e'\in E_+$ with $o(e')=t(e)=v$, so that $t(e')=yv$ and $o(e)=x^{-1}v$ for (unique) $y,x\in P$. Then
   $\Omega(e') \subset \Omega(e)$ if and only if $y\notin \lambda(x)$.
   This is because $\Omega(e') \subset \Omega(e)$ if and only if
   the edges $e$ and $e'$ lie as shown in Figure \ref{hexagon}, in some apartment.
\refstepcounter{picture}
\begin{figure}[htbp]\label{hexagon}
\centerline{
\beginpicture
\setcoordinatesystem units <0.6cm,1cm> 
\setplotarea x from -5 to 5, y from 1 to 3         
\setlinear
\plot -1 1  1 1   2 2  1 3  -1 3  0 2  1 3 /
\plot -1 1  0 2  1 1  /
\plot  -1 1  -2 2  -1 3 /
\plot -1.5 1  1.5 1  / 
\plot  -2.2 2  2.5 2 /  
\plot -1.5 3  1.5 3  /   
\plot -1.2 0.8  1.2 3.2  /
\plot  1.2 0.8  -1.2 3.2 /
\plot 2.2 1.8  0.8 3.2  /
\plot -2.2 1.8  -0.8 3.2  /
\plot -2.2 2.2  -0.8 0.8  /
\plot 2.2 2.2  0.8 0.8  /
\put{$_e$}    at  -0.9 2.2
\put{$_{e'}$}    at  1.1 2.3
\arrow <10pt> [.2, .67] from  -1.3 2 to -0.8 2
\arrow <10pt> [.2, .67] from  0.3 2 to 1.2 2
\endpicture
}
\caption{}
\end{figure}

\noindent Equations (\ref{boundary A}) and (\ref{boundary C}) imply the following relations in
$C(\cl P,\bb Z)_\Gamma$.

\begin{subequations}\label{coinv C}
\begin{eqnarray}
\sum_{y\notin \lambda(x)}[y] &=& [x], \quad\, x\in P; \label{A}\\
\sum_{x\in P}[x] &=& [\id].\label{C}
\end{eqnarray}
\end{subequations}

These should be compared with relations (\ref{a}),(\ref{c}), respectively.
Now we seek an analogue of (\ref{b}). The following fact is the key.

\begin{lemma}\label{ronan}
{\rm \cite[Lemma 9.4]{ron}}
Given any chamber $c$ and any sector $S$ in $\Delta$, there exists a sector $S_1\subset S$ such that $S_1$ and $c$ lie in a common apartment.
\end{lemma}

\begin{figure}[htbp]
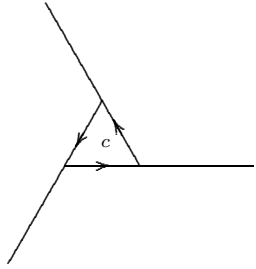
\label{chamberray}
\hfil
\centerline{
\beginpicture
\setcoordinatesystem units <0.5cm,0.866cm> 
\setplotarea x from -4 to 4, y from 0 to 4 
\setlinear \plot  -1.5 3.5  1 1 /
\put {$_{c}$} [] at 0.1 1.35
\putrule from -1 1 to 4 1
\setlinear \plot  -2.5 -0.5    0 2  /
\arrow <5pt> [.2, .67] from   0 1   to 0.2 1
\arrow <5pt> [.2, .67] from  -0.5 1.5    to -0.7 1.3
\arrow <5pt> [.2, .67] from  0.5 1.5     to 0.3 1.7
\endpicture
}
\hfil
\caption{Relative positions of a chamber and representative rays}
\end{figure}

If $c$ is a chamber of $\Delta$ and if $\xi\in\cl P$, then $\xi$ has a representative ray that lies relative to $c$ in one of the three positions in Figure \ref{chamberray}, in some apartment containing them both. This is because, by Lemma \ref{ronan}, we can choose a ray $[x,\xi)$ such that $c$ and
$[x,\xi)$ lie in a common apartment. Now choose an appropriate ray parallel to $[x,\xi)$.
\refstepcounter{picture}

\noindent The next lemma follows immediately. Equation (\ref{B}) is the
desired analogue of (\ref{b}).
\begin{lemma}\label{3}
  If $e_0, e_1,e_2\in E_+$ are the edges of a chamber $c \in \Delta$ then
  $$\Omega(e_0)\sqcup\Omega(e_1)\sqcup\Omega(e_2) = \cl P.$$
  Consequently, if $\cl T$ is a triangle presentation and $(x,y,z) \in \cl T$ then
 \begin{equation}\label{B}
      [x]+[y]+[z] = [\id].
 \end{equation}
\end{lemma}

Now we specify generators for $C(\cl P,\bb Z)_\Gamma$.

\begin{lemma}\label{fg}
The group $C(\cl P,\bb Z)_\Gamma$ is finitely generated, with generating set $\{[x]:x\in\cl P\}$.
\end{lemma}

\begin{proof}
Every clopen set $V$ in $\cl P$ may be expressed as a finite disjoint union of sets of the form $\Omega(e)$, $e\in E_+$. Any function $f\in C(\cl P,\bb Z)$ is bounded, by compactness of $\cl P$, and so takes finitely many values $n_i\in\bb Z$. Therefore $f$ may be expressed as a finite sum $f=\sum_j n_j\chi_{e_j}$, with
$e_j\in E_+$. The result follows.
\end{proof}

\begin{proposition}\label{hom}
  There is a homomorphism $\theta$ from $A_{\cl T}$ onto $C(\cl P,\bb Z)_\Gamma$ defined by
  $\theta(x)=[x]$, for $x\in P$ and $\theta(\varepsilon)=[\id]$.
\end{proposition}

\begin{proof}
  Equations (\ref{A}),(\ref{C}),(\ref{B}) show that $\theta$ preserves the relations (\ref{a}),(\ref{c}),(\ref{b}), respectively.
Therefore $\theta$ extends to a homomorphism. Surjectivity is a consequence of Lemma \ref{fg}.
\end{proof}

\section{The main results}

This section collects the main consequences.

\begin{corollary}\label{C1} Let $\Gamma$ be an $\widetilde A_2$ group acting on an $\widetilde A_2$ building $\Delta$ of order $q$ with minimal boundary $\cl P$.
Then $C(\cl P,\bb Z)_{\Gamma}$ is a finite group and the class $[\id]$  in $C(\cl P,\bb Z)_{\Gamma}$ has order bounded by $q^2-1$. If $\Gamma$ has an S-invariant triangle presentation then $(q-1)[\id]=0$.
\end{corollary}

\begin{proof}
  This follows immediately from Proposition \ref{finite}, Proposition \ref{hom}, Lemma \ref{q2} and
   Corollary \ref{S-invariant}.
\end{proof}

The final statement of Corollary \ref{C1} applies, in particular, to the
groups of Tits type (Section \ref{appendix}).

\subsection{K-theory}\label{KT}

 The {\it Furstenberg boundary} $\Omega$ of $\Delta$ is the set of chambers of $\Delta^\infty$, endowed with a compact totally disconnected topology in which basic open sets have the form
$$
\Omega(v) = \left \{ \omega \in \Omega : [x,\omega) \ \text{ contains } v \right \}
$$
where $v$ is a vertex of $\Delta$ and $[x,\omega)$ is the representative sector for $\omega$ with base vertex $x$ \cite[Section 2]{cms}.
If $\Delta$ is the Bruhat-Tits building of $G=\PGL_3(\bb K)$ then $\Omega$ is isomorphic as a topological $G$-space to the space of maximal flags $(0)<V_1<V_2<\bb K^3$.
The mapping which sends each sector to its wall of type $0$
induces a natural surjection $\Omega\to\cl P$, under which $\cl P$
has the quotient topology. Since this surjection is equivariant, there is an induced epimorphism $C(\cl P ,\bb Z)_{\Gamma}\to C(\Omega,\bb Z)_{\Gamma}$.

The topological action of an $\widetilde A_2$ group $\Gamma$ on the maximal boundary is encoded in the full crossed product $C^*$-algebra $\fk A_{\Gamma}=C(\Omega)\rtimes\Gamma$, which is studied in \cite{rs1,rs2,rs3}.
The natural embedding  $C(\Omega) \to \fk A_{\Gamma}$ induces a homomorphism
\begin{equation}\label{psi}
  \psi:C(\Omega,\bb Z)_\Gamma\to K_0(\fk A_{\Gamma})
\end{equation}
and $\psi([\id])=[\id]_{K_0}$, the class of $\id$ in the $K_0$-group of $\fk A_{\Gamma}$.
The article \cite{em} estimates the order of $[\id]_{K_0}$ for various boundary $C^*$-algebras, and contains an extensive bibliography.

In \cite{rs3}, T. Steger and the second author performed extensive computations which determined the order of $[\id]_{K_0}$ for  many $\widetilde A_2$ groups with $q\le 13$. The computations were done for all the $\widetilde A_2$ groups in the cases $q=2,3$ and for several representative groups for each of the other values of $q\le 13$.
If $q=2$ there are precisely eight $\widetilde A_2$ groups $\Gamma$, all of which embed as lattices in a linear group $\PGL (3,\bb K)$ where $\bb K= \bb F_2((X))$ or $\bb K=\bb Q_2$. If $q=3$ there are 89 possible $\widetilde A_2$ groups, of which 65 do not embed naturally in linear groups.
The experimental evidence suggested that for boundary crossed product algebras associated with $\widetilde A_2$ groups it is always true that
$[\id]_{K_0}$ has order $q-1$ for $q \not\equiv 1 \pmod {3}$ and has order $(q-1)/3$ for $q  \equiv 1 \pmod {3}$. It is striking that the order of $[\id]$ appears to depend only on the parameter $q$.
It is shown in \cite{rs3} that the order of $[\id]_{K_0}$ is bounded above by $q^2-1$ and below by $\frac{q-1}{(q-1,3)}$.

\begin{corollary}\label{C*}
Let $\cl T$ be an S-invariant triangle presentation and $\Gamma=\Gamma_{\cl T}$.
Then the class $[\id]_{K_0}$  in $K_0(\fk A_{\Gamma})$ has order bounded by $q-1$.
If $q \not\equiv 1 \pmod {3}$ then the order of $[\id]_{K_0}$ is precisely $q-1$.
\end{corollary}

\begin{proof}
This follows directly from Corollary \ref{C1}, since $\psi([\id])=[\id]_{K_0}$.
\end{proof}

Theorem \ref{main} is an immediate consequence of this Corollary.

\begin{remark}
It also follows from the computations in Remark \ref{T0'} that the order of $[\id]_{K_0}$ is $\frac{q-1}{3}$, for all groups of Tits type with $q \equiv 1 \pmod {3}$ and $q\le 31$.
\end{remark}

\begin{remark}
  For each $\widetilde A_2$ group $\Gamma$, the algebra $\fk A_{\Gamma}$ has the structure of a higher rank Cuntz-Krieger algebra \cite[theorem 7.7]{rs2}. These algebras are classified up to isomorphism by their two $K$-groups, together with the class $[\id]_{K_0}$, \cite[Remark 6.5]{rs2}. It was proved in \cite[Theorem 2.1]{rs3} that
  \begin{equation}\label{Kremark}
  K_0(\fk A_{\Gamma})=K_1(\fk A_{\Gamma})=\bb Z^{2r}\oplus T,
  \end{equation}
  where $r\ge 0$ and $T$ is a finite group. The computations done in \cite{rs3}
  give rise to some striking observations. For example, there are precisely three torsion-free $\widetilde A_2$ subgroups of $\PGL_3(\bb Q_2)$ and these three groups are distinguished from each other by $K_0(\fk A_{\Gamma})$.

  There is also strong evidence that, for any torsion free
  $\widetilde A_2$ group $\Gamma$,
  the integer $r$ in (\ref{Kremark}) is equal to the rank of $H_2(\Gamma, \bb Z)$.
  A typical example is provided by the group $\Gamma=\Gamma_{\cl T_0}$,
  associated with the triangle presentation of Tits type defined in Section \ref{appendix},
  with $q=13$. In that case,
    $$K_0(\cl A_\Gamma)=\bb Z^{1342}\oplus(\bb Z/4\bb Z)\oplus(\bb Z/3\bb Z)^{6}
    \oplus(\bb Z/13\bb Z)^{6}.$$
    As expected, the class $[\id]_{K_0}$ has order $4=\frac{q-1}{3}$ and
    $H_2(\Gamma, \bb Z)=\bb Z^{671}$.

\end{remark}


\bigskip

\begin{thebibliography}{}

\bibitem{ab} P. Abramenko and K. Brown, {\it Buildings. Theory and applications}, Graduate Texts in Mathematics, 248. Springer, New York, 2008.

\bibitem{bhv} B. Bekka, P. de la Harpe and A. Valette, {\it Kazhdan's Property (T)}, Cambridge University Press, Cambridge, 2008.

\bibitem{cms}  D. I. Cartwright, W. M{\l}otkowski and T. Steger, Property (T) and
$\widetilde A_2$ groups, {\it Ann.\ Inst.\ Fourier} {\bf 44} (1993), 213--248.

\bibitem{cmsz} D. I. Cartwright, A. M. Mantero, T. Steger and
A. Zappa, Groups acting simply transitively on the vertices of a
building of type~$\widetilde A_2$, I and II,\ {\it Geom.\ Ded.}  {\bf 47}
(1993), 143--166 and 167--223.

\bibitem{em} H. Emerson and R. Meyer,
Euler characteristics and Gysin sequences for group actions on boundaries,
{\it Math. Ann.} {\bf 334} (2006), 853--904.

\bibitem{rs1} G. Robertson and T. Steger, $C^*$-algebras arising from group actions on the boundary of a triangle building, {\it Proc.\ London Math.\ Soc.}  {\bf 72} (1996), 613--637.

\bibitem{rs2} G. Robertson and T. Steger, Affine buildings, tiling systems and higher rank Cuntz-Krieger algebras, {\it J. reine angew. Math.} {\bf 513} (1999), 115--144.

\bibitem{rs3} G. Robertson and T. Steger, Asymptotic K-theory for groups acting on $\tilde A_2$ buildings, {\it Canad. J. Math.} {\bf 53} (2001), 809--833.

\bibitem{ron} M. Ronan, {\it Lectures on Buildings}, University of Chicago Press, 2009.

\bibitem{sin} J. Singer, A theorem in finite projective geometry and some applications,
{Trans. Amer. Math. Soc.} {\bf 43} (1938), 377--385.

\bibitem{Wei} R.~Weiss, {\it The Structure of Affine Buildings}, Annals of Mathematics Studies, Vol.~168,
Princeton, 2009.

\bibitem{Zuk} A.~Zuk, La propri\'et\'e (T) de Kazhdan pour les groupes agissant sur les poly\`edres, {\it C. R. Acad. Sci. Paris S\'er. I Math.}  {\bf 323}  (1996), 453--458.

\end{thebibliography}
\end{document}